\documentclass{article}
\usepackage[T2A]{fontenc}
\usepackage[cp1251]{inputenc}
\usepackage[english,russian]{babel}
\usepackage[tbtags]{amsmath}
\usepackage{amsfonts,amssymb,mathrsfs,amscd}
\usepackage{cmap,amsmath,amsthm,hyperref}
\allowdisplaybreaks
\makeatletter
\gdef\No{{\select@language{russian}\textnumero}}
\makeatother

\theoremstyle{plain}
\newtheorem{theorem}{Теорема}
\newtheorem{lemma}{Лемма}[section]

\theoremstyle{definition}

\def\const{\mathrm{const}}

\newcommand{\itemaa}{\smallskip${{\bigstar}}$\hspace{0.1cm}}
\newcounter{inum}
 
\newcounter{meq}
\def\mmeq#1{\refstepcounter{meq}\begin{equation*} #1\eqno{(\hbox{S}\themeq)}\end{equation*}}

\def\seqref#1{(S\ref{#1})}

\def\eqa#1{\begin{equation}\begin{aligned}#1\end{aligned}\end{equation}}
\def\eq#1{\begin{equation}#1\end{equation}}
\def\eqs#1{\begin{eqnarray}#1\end{eqnarray}}

\def\seq#1{\begin{equation*}#1\end{equation*}}
\def\seqs#1{\begin{equation*}\begin{aligned}#1\end{aligned}\end{equation*}}

\def\phi{\varphi}

\def\be{\begin{equation}}
\def\ee{\end{equation}}
\def\ba{\begin{aligned}} % аналог array, но более удобный
\def\ea{\end{aligned}}

\def\const{\mbox{const}}

\newcounter{theo}

\newcounter{lem}

\newcounter{prop}
\newcounter{rem}

\newcounter{defi}

\newcounter{examp}

\def\proof{\noindent {\bf Доказательство.} }
\def\pfrac#1#2{\frac{\partial #1}{\partial #2}}

\begin{document}
\title{Гиперболические уравнения с симметриями пятого порядка.}

\author{Р.\,Н.~Гарифуллин}%\thanks{\rm Исследования в Секциях 1-2, 4 выполнено за счет гранта Российского научного фонда № 21-11-00006, https://rscf.ru/project/21-11-00006/.}
%\address{Институт математики с вычислительным центром Уфимского федерального исследовательского центра РАН, Уфа, Россия.}
%\email{rustem@matem.anrb.ru}
%второй автор
%\author[I.\,I.~Ivanov]{И.\,И.~Иванов}
%\address{}
%\email{}

\date{}
%\udk{517.929}

\maketitle
%\begin{fulltext}

\begin{abstract}
 В работе проводится классификация уравнений гиперболического типа. Исследуется класс уравнений вида
\seq{\frac{\partial^2 u}{\partial x\partial y}=F\left(\frac{\partial u}{\partial x},\frac{\partial u}{\partial y},u\right),} здесь $u(x,y)$ -- неизвестная функция, $x,y$ независимые переменные. 
Классификации основывается на требовании существования высших симметрий пятого порядка. В результате получен список 4 уравнений с требуемыми условиями.

\end{abstract}

%\begin{

%\markright{Классификация полудискретных уравнений гиперболического типа.}% Случай симметрий пятого порядка.}

%\footnotetext[0]{Исследования в Секциях 1-2, 4 выполнено за счет гранта Российского научного фонда №\,21-11-00006, https://rscf.ru/project/21-11-00006/.}

В этой работе исследуются  уравнения гиперболического типа\eq{u_{xy}=F\left(u_x,u_y,u\right)\label{eq},} где неизвестная функция $u(x,y)$ зависит от  двух переменных $x$ и $y$. Среди интегрируемых уравнений вида \eqref{eq} выделяются три следующих уравнения \cite{zs79}:
\eqs{&&u_{xy}=e^u,\label{darb}\\&&u_{xy}=e^u-e^{-u}\label{sin},\\ &&u_{xy}=e^{u}+e^{-2u}\label{tz}.} Первое их этих уравнений является интегрируемым по Дарбу, и задача поиска полного списка таких уравнений рассматривалась в работе \cite{zs01}, см. также \cite{vz23,s25}. Второе уравнение \eqref{sin} -- это известное уравнение синус-Гордона, известно что это уравнение допускает высшие симметрии третьего порядка. Классификация уравнений \eqref{eq} с симметриями третьего порядка проведена в статье \cite{ms11}. Уравнение \eqref{tz} -- это уранения Цицейки, оно не допускает высших симметрий третьего порядка, но обладает высшими симметриями пятого порядка: \eqa{&u_t=u_{xxxxx}+5(u_{xx}-u_x^2)u_{xxx}-5u_xu_{xx}^2+u_x^5,\label{sxy} \\ &u_{\tau}=u_{yyyyy}+5(u_{yy}-u_y^2)u_{yyy}-5u_yu_{yy}^2+u_y^5.} Насколько известно автору классификации уравнений \eqref{eq} с симметриями пятого порядка не проводились. Ранее были были известны только примеры уравнений такого вида \cite{zs01,bzp02}.

Из результатов работы \cite{z94} можно показать (см. также \cite{ms11}), что если уравнение \eqref{eq} допускает высшую симметрию вида $$u_t=H(u_{5},u_{4},u_{3},u_2,u_1,u),$$ где обозначено $u_j=\pfrac{^ju}{x^j},$ то функция $G$ имеет вид уравнению \eq{H=u_{5}+G(u_{4},u_{3},u_2,u_1,u). \label{c_sep}} Во всех известных примерах высшие симетрии сами по себе являются интегрируемыми уравнениями, см. \cite{zs01,ms11} для непрерырвных уравнений, \cite{g23,g25} для полудискретных уравнений, \cite{x09,ly11,gy12} для полностью дискретных уравнений. Поэтому естественно предполагать, что уравнения \eqref{c_sep} должно быть интегрируемым. Уравнения \eqref{c_sep} линейные по старшей производной называются уравнениями с постоянной сепарантой и полный список интегрируемых уравнений такого вида приведен в работах \cite{dss85,ms12}, во второй из которых также подробно изложена история вопроса. 

Список таких уравнений пятого порядка имеет вид:
\begin{align}
&u_t=u_{5}+5 u u_{3}+5 u_1 u_2+5 u^2 u_1, \label{tst1} \\
\label{tst2}
&u_t=u_{5}+5 u u_{3}+{25 \over 2} u_1 u_2+5 u^2 u_1,  \\
\label{tst3}
&u_t=u_{5}+5 u_1 u_{3}+{5 \over 3} u_1^3,  \\ 
\label{tst4}
&u_t=u_{5}+5 u_1 u_{3}+{15\over 4}u_2^2 + {5 \over 3} u_1^3, \\
\label{tst5}
&u_t=u_{5}+5 (u_1-u^2) u_{3}+5 u_2^2-20 u u_1 u_2-5 u_1^3+5 u^4 u_1, \\
\label{tst6}
&u_t=u_{5}+5 (u_2-u_1^2) u_{3}-5 u_1 u_2^2+u_1^5, \\
\label{tst7}
&\begin{aligned}
u_t&=u_{5}+5 (u_2-u_1^2-\lambda_1^2 e^{2u}-\lambda_2^2 e^{-4u}) \, u_{3}-5 u_1 u_2^2\\&-15 (\lambda_1^2 e^{2u}-4 \lambda_2^2 e^{-4u})\, u_1 u_2 
+u_1^5-90 \lambda_2^2 e^{-4u}\, u_1^3\\&+5(\lambda_1^2 e^{2u}+\lambda_2^2 e^{-4u})^2\, u_1,
\end{aligned}\\
\label{tst8}
&\begin{aligned}
u_t&=u_{5}+5 (u_2-u_1^2-\lambda_1^2 e^{2u}+\lambda_2 e^{-u}) \, u_{3}-5 u_1 u_2^2-15 \lambda_1^2 e^{2u} \, u_1 u_2 \\
&+u_1^5+5(\lambda_1^2 e^{2u}-\lambda_2 e^{-u})^2 \, u_1, \quad \lambda_2\ne 0, 
\end{aligned}\\
\label{tst9}
&\begin{aligned}
u_t&=u_{5}-5\frac{u_2 u_{4} }{ u_1}+ 5\frac{u_2^2 u_{3}}{ u_1^2}+5\left(\frac{ \mu_1}{ u_1}+\mu_2 u_1^2\right)u_{3}-5\left(\frac{ \mu_1}{ u_1^2}+\mu_2 u_1\right)u_2^2
\\
&-5 \frac{\mu_1^2}{ u_1}+ 5 \mu_1\mu_2 u_1^2 +\mu_2^2 u_1^5,
\end{aligned}%\\
\end{align}
\begin{align}
\label{tst10}
&\begin{aligned}
u_t&=u_{5}-5\frac{u_2 u_{4}}{u_1}-\frac{15}{4 }\,\frac{u_{3}^2}{u_1}+\frac{ 65}{4}\,\frac{u_2^2 u_{3}}{ u_1^2} +5\left(\frac{\mu_1}{u_1}+\mu_2 u_1^2\right)\, u_{3}
-\frac{135}{16 }\frac{u_2^4}{u_1^3} 
\\
&-5\left(\frac{7 \mu_1}{4 u_1^2}-\frac{\mu_2 u_1}{2}\right) u_2^2-5 \frac{\mu_1^2}{u_1}+ 5 \mu_1\mu_2 u_1^2 +\mu_2^2 u_1^5,
\end{aligned}\\
\label{tst11}
&\begin{aligned}
u_t&= u_{5}-\frac{ 5 }{2}\,\frac{ u_2u_{4} }{ u_1}-\frac{5}{4}\,\frac{u_{3}^2}{ u_1}+5\frac{u_2^2 u_{3}}{u_1^2} +\frac{5\, u_2 u_{3}}{2 \sqrt{u_1}}
-5 (u_1-2 \mu u_1^{1/2}+\mu^2) \, u_{3} \\&-\frac{ 35}{16}\,\frac{u_2^4}{ u_1^3}  
  -\frac{5}{3}\,\frac{u_2^3}{u_1^{3/2}} +5\Big (\frac{ 3 \mu^2}{ 4 u_1} -\frac{\mu }{\sqrt{u_1}}+\frac{1}{4}\Big)\, u_2^2
+ \frac{ 5}{ 3}\, u_1^3  
 - 8 \mu u_1^{5/2}\\&+15 \mu^2 u_1^2-\frac{40}{3}\,\mu^3 u_1^{3/2},
\end{aligned}\\
&\begin{aligned}\label{tst12m}
u_t&=u_{{{ 5}}}+\frac52\,{\frac { f-u_{1}}{{f}^{2}}}\,u_{{2}}u_{{{ 4}}}+\frac54\,{\frac {2\,f -u_{1}}{{f}^{2}}}\,u_3^{2} +5\,\mu\, ( u_{1}+f ) ^{2}u_{{{3}}}\\
&+\frac54\,{\frac {4\,{u_{1}}^{2}-8\,u_{1}f+{f}^{2}}{{f}^{4}}}\,u_2^{2}u_{{{ 3}}} +{\frac {5}{16}}\,{\frac {2-9\,u_{1}^{3}+18\,u_{1}^{2}f}{{f}^{6}}}\,u_2^{4}\\
&+\frac {5\mu}{4}\,{\frac {( 4\,f-3\,u_{1} )( u_{1}+f )^{2}}{{f}^{2}}}\,u_2^{2}+{\mu}^{2} ( u_{1}+f ) ^{2}\big( 2\,f ( u_{1}+f ) ^{2}-1\big),
\end{aligned}\\
&\begin{aligned}\label{tst13m}
u_t&=u_{{{ 5}}}+\frac52\,\frac {f- u_{1}}{f^{2}}\, u_2u_4+\frac54\,\frac {2\,f- u_1}{f^2}\, u_3^{2}-5\,\omega\, ( {f}^{2}+u_1^{2} ) u_3 \\
&+\frac54\,{\frac {4\,u_1^{2}-8\,u_{1}f+{f}^{2}}{{f}^{4}}}\,u_2^{2}u_3+{\frac {5}{16}}\,\frac {2-9\,u_1^{3}+18\,u_1^{2}f}{f^6}\,u_2^{4}
\\
&+\frac54\,\omega\,\frac {5\,u_1^{3} -2\,u_1^{2}f-11\,u_{1}{f}^{2}-2}{f^2}\, u_2^{2}-\frac52\,{\omega'}\, ( u_1^{2}-2\,u_{1}f+5\,{f}^{2} )u_{1}u_2 \\
&+5\,{\omega}^{2}u_{1}{f}^{2} ( 3\,u_{1}+f ) ( f-u_{1}),
\end{aligned}
\end{align}
\begin{align}
&\begin{aligned}\label{tst14m}
u_{{t}}&=u_5+\frac52\,{\frac {f-u_{1}}{{f}^{2}}}\,u_2u_4+\frac54\,{\frac {2\,f -u_{1}}{{f}^{2}}}\,u_3^{2}
+\frac54\,{\frac {4\,u_1^{2}-8\,u_{1}f+{f}^{2}}{{f}^{4}}}\,u_2^{2}u_3
\\
&+{\frac {5}{16}}\,{\frac {2-9\,{u_{1}}^{3}+18\,{u_{1}}^{2}f}{{f}^{6}}}\,u_2^{4}+5\,\omega\,{\frac {2\,{u_{1}}^{3}+{u_{1}}^{2}f-2\,u_{1}{f}^{2}+1}{{f}^{2}}}\,u_2^{2}
\\
& -10\,\omega\,u_{{{3}}} ( 3\,u_{1}f+2\,u_1^{2}+2\,{f}^{2} )-10\,{\omega'} ( 2\,{f}^{2}+u_{1}f+{u_{1}}^{2} )\,u_{1}u_{{2}}\\
& +20\,{\omega}^{2}u_{1} ( {u_{1}}^{3}-1 )( u_{1}+2\,f ),
\end{aligned}
\end{align}
\begin{align}
\label{tst15m}
&\begin{aligned}
u_t&=u_5+\frac52\,\frac {f-u_1}{f^2}\,u_2u_4+\frac54\,{\frac {2\,f-u_{1}}{{f}^{2}}}\,u_3^{2}-5\,c\frac {{f}^2+u_1^2}{{\omega}^{2}}\,u_3
\\
&+\frac54\,{\frac {4\,u_1^{2}-8\,u_{1}f+{f}^{2}}{{f}^{4}}}\,u_2^{2}u_3 +\frac {5}{16}\,\frac {2-9\,{u_{1}}^{3}+18\,u_1^{2}f }{{f}^{6}}\,u_2^{4}
\\
&-10\,\omega\, ( 3\,u_{1}f+2\,u_1^{2}+2\,{f}^{2} )\,u_3-\frac54\,c\,\frac { 11\,u_{1}{f}^{2}+2\,u_1^{2}f +2-5\,u_1^{3} }{{\omega}^{2}{f}^2}\,u_2^{2}
\\
&+5\,\omega\,\frac {2\,u_1^{3}+u_1^{2}f-2\,u_{1}{f}^{2}+1}{f^2}\,u_2^{2}+5\,c\,{\omega'}\,\frac {u_1^{2}+5\,{f}^{2}-2\,u_{1}f }{{\omega}^3}\,u_{1}u_2
\\
&-10\,{\omega'}\,(2\,{f}^{2}+u_{1}f+{u_{1}}^{2})\,u_1u_2 +20\,{\omega}^{2}u_{1} ( u_1^{3}-1 )  ( u_{1}+2\,f ) 
\\
&+40\,{\frac {c\,u_{1}{f}^{3} ( 2\,u_{1}+f ) }{\omega}}+5\,{\frac {{c}^{2}u_{1}{f}^{2} ( 3\,u_{1}+f )( f-u_{1} ) }{{\omega}^{4}}},\ \ c\ne0.
\end{aligned}
\end{align}
{\it Здесь} $\lambda_1, \lambda_2,\mu, \mu_1, \mu_2$ {\it и} $c$ --- {\it параметры},
%the expression $ \{\lambda_1,\lambda_2\}$ is a vector. 
{\it функция}  $f(u_1)$ {\it является решением алгебраического уравнения}
\begin{equation}\label{alg1}
(f+u_1)^2(2f-u_1)+1=0,
\end{equation}
{\it а} $\omega(u)$ --- {\it это любое непостоянное решение дифференциального уравнения}
\begin{equation}\label{Waier}
\omega '^2=4\, \omega^3+c. \qquad \square
\end{equation}

Приведенный список уравнений отличается от списка работы \cite{ms12} переобозначением константы $\lambda_1$ в уравнение \eqref{tst7}. В списке уравнений \eqref{tst1} -- \eqref{tst15m} нижний числовой индекс обозначает порядок производной по $x$.

Высшие симметрии во втором направлении являются интегреруемыми уравнениями вида:
\eq{u_{\tau}=u_{yyyyy}+G(u_{yyyy},u_{yyy},u_{yy},u_{y},u),\label{y_sym}} и также должны принадлежать списку \eqref{tst1} -- \eqref{tst15m} при замене $x\leftrightarrow y$ с точностью до точечных замен.

%В работе ищутся только автономные по дискретной переменной $n$ уравнения вида \eqref{seq} и соответственно используются только автомномные высшие симметрии (\ref{tst1}--\ref{tst15m}). Это, с одной стороны, является упрощением задачи, но, с другой стороны, не известно неавтономных уравнений вида \eqref{seq} совместных с дискретной высшей симметрией (возможно неавтономной) вида \eqref{d_sym}. Отметим, что полученные ответы содержат дополнительные произвольные константы, вместо этих констант можно брать функции зависящие от $n$, при этом останется совместность с непрерывными высшими симметриями. 

%В полностью диcкретном случае, сущуствуют автономные уравнения, у которых одна из высший симметрий имеет порядок 3 и является неавтомной, а высшая симметрия в другом направление имеет больший порядок, см. \cite{gy12,gmy14,gy19}. Подобные уравнения в этой работе не исследуются. Вопрос об их существование остается открытым.

\section{Метод исследования.}

Из требования совместности уравнений \eqref{eq} и \eqref{c_sep} получаем определяющее уравнение
\eq{D_xD_y (u_5+G)=\pfrac{F}{u_{x}}D_x(u_5+G)+\pfrac{F}{u_{y}}D_y(u_5+G)+\pfrac{F}{u}(u_5+G). \label{lineq}}
Здесь $D_x, D_y$ операторы полных производных по переменных $x,y$ соответственно. Исключая из \eqref{lineq} смешанные производные в силу уравнения \eqref{eq} получаем соотношение, которое должно тождественно выполняться по переменным $u,u_{1},u_{2},u_3,u_4,u_5,u_y$. 
При $u_5$ получаем уравнение \eq{D_y\pfrac{G}{u_4}+5D_x\pfrac{F}{u_x}=0.\label{eq5}} 

Исследование уравнение \eqref{eq5} приводит к лемме
\begin{lemma}Пусть  $\pfrac{G}{u_4}=5u_2g(u_1)$, %где $g(u_1)$ представляется одним из следующих видов: $$g(u_1)=0,\quad g(u_1)=-\frac{1}{u_1},\quad g(u_1)=-\frac12\frac{1}{u_1},\quad g(u_1)=\frac12\,{\frac {f-u_{1}}{{f}^{2}}}. $$Тогда  
тогда решение уравнения \eqref{eq5} имеет вид: \eq{F=F_1(u_y,u)f_1(u_1)+C_2f_2(u_1),\label{ras}} где $F_1$ произвольная функция, $f_1,f_2$ линейно независимые решения уравнения $w''+gw'+g'w=0$, $C_2$ константа.
\end{lemma}
\proof При условиях Леммы уравнение \eqref{eq5} можно разделить на 2 уравнения по переменной $u_2$:
\eqs{&\pfrac{^2F}{u_1^2}+g(u_1)\pfrac{F}{u_1}+g'(u_1)F=0,\label{eqF1}\\
&u_1\left(\pfrac{^2F}{u_1\partial u}+g(u_1)\pfrac{F}{u}\right)+F\left(\pfrac{^2F}{u_1\partial u_y}+g(u_1)\pfrac{F}{u_y}\right)=0.\label{eqF2}}
Уравнение \eqref{eqF1} можно проинтегровать 1 раз:
\seq{\pfrac{F}{u_1}+g(u_1)F=F_2(u_y,u)} и решить:
\seq{F=\left(F_2(u_y,u)\int \exp\left(\int g(u_1)du_1\right)du_1+F_1(u_y,u)\right)\exp\left(-\int g(u_1)du_1\right).} Подстановка этого представления  в уравнение \eqref{eqF2} приводит к:
\eqa{F_2\pfrac{F_2}{u_y}\,\exp\left(-\int g(u_1)du_1\right)\int \exp\left(\int g(u_1)du_1\right)du_1\,\\+F_1\pfrac{F_2}{u_y}\exp\left(-\int g(u_1)du_1\right) +u_1\pfrac{F_2}{u_0}=0.\label{eq3s}} Вронскин трех функций $$\exp\left(-\int g(u_1)du_1\right)\int \exp\left(\int g(u_1)du_1\right)du_1,  \exp\left(-\int g(u_1)du_1\right),u_1$$ показывает, что они зависимы только при $g(u_1)=\frac{a}{u_1}$, поэтому во всех остальных случаях получаем $\pfrac{F_2}{u_y}=\pfrac{F_2}{u_0}=0$, следовательно $F_2=C_2$. При $g(u_1)=\frac a{u_1}$ уравнение \eqref{eq3s} принимает вид:
\seq{\left(\frac{1}{a+1}F_2\pfrac{F_2}{u_y}+\pfrac{F_2}{u_0}\right)u_1+u_1^{-a}F_1\pfrac{F_2}{u_y}=0, \hbox{ при }a\neq-1, a\neq 0,}   
\seq{\left(F_1\pfrac{F_2}{u_y}+\pfrac{F_2}{u_0}\right)u_1+u_1\ln u_1 F_2\pfrac{F_2}{u_y}=0,\hbox{ при }a=-1,}   \seq{\left(F_2\pfrac{F_2}{u_y}+\pfrac{F_2}{u_0}\right)u_1+F_1\pfrac{F_2}{u_y}=0,\hbox{ при }a=0.} Во всех трех случаях видим, что либо $F_2=\const$, либо $F_1=0$ и утверждение Леммы верно с точностью до переобозначения $F_2\to F_1$.\qed

Замечание. Можно заметить, что все уравнения \eqref{tst1} -- \eqref{tst15m} удовлетворяют условию Леммы. Для них функции $g(u_1)$ и $F$ имеют один из видов:
\seqs{&g(u_1)=0,\quad F=F_1(u_y,u_0)+C_2u_1 \hbox{ или } F=F_1(u_y,u_0)u_1,\\
&g(u_1)=-\frac{1}{u_1},\quad F=F_1(u_y,u_0)u_1+C_2u_1\ln u_1,\\ &g(u_1)=-\frac12\frac{1}{u_1},\quad F=F_1(u_y,u_0)\sqrt{u_1}+C_2u_1 \hbox{ или } F=F_1(u_y,u_0)u_1,\\ &g(u_1)=\frac12\,{\frac {f-u_{1}}{{f}^{2}}},\quad F=F_1(u_y,u_0)f(u_1)+C_2\ln(f(u_1)+u_1). }

Дальнейшее исследование уравнения \eqref{lineq} при конкретных функция $G$ из \eqref{tst1} -- \eqref{tst15m} позволяет найти функции $F_1$ и показывают, что $C_2=0.$

Подобный подход к систематическому поиску интерируемых уравнений гиперболического типа известными симметриями применялся ранее в непрерывном случае в работе \cite{ms11} для скалярных уравнений с симметриями третьего порядка, в статье \cite{ms14} для векторных уравнений, в полудискретном случае в работах \cite{g23,g25} для симметрий третьего и пятого порядка и в работах \cite{gy15, ggy19} для полностью дискретных уравнений.

В ходе подобной классификации также возникают Дарбу интегрируемые уравнения вида \eqref{eq}, такие уравнения наряду с высшими симметриями обладают двумя нетривиальными интегралами (см. \cite{zs01}) вдоль характеристических направлений:
\seq{W_1(u,u_{1},u_{2},\ldots,u_{k}),\quad W_2(u,u_{y},u_{yy},..\ldots, \pfrac{^k}{y^k}u),} такими что \seq{D_y W_1=0,\quad D_x W_2=0.}
В ходе проведенной классификации найденные уравнения такого типа оказались ранее известными \cite{zs01}, поэтому здесь они не были продемонстрированы.

\section{Результаты классификации.}
В данной секции приводятся найденные гиперболические уравнения. Они сгруппированы по виду высшей симметрии в $x$ направлении. Верно следующее утверждение:
%\begin{itemize}
\begin{theorem} Если невырожденное нелинейное автономное уравнение \eqref{eq} допускает высшую симметрию в $x$ направлении в виде одного из уравнений (\ref{tst1} -- \ref{tst15m}) и не является Дарбу интегрируемым, то оно имеет вид \seqref{tzm} -- \seqref{tz01}.

\end{theorem}

{\bf Схема доказательства. } Правые части уравнений (\ref{tst1}--\ref{tst15m}) брались в качестве функции $V_n$ в определеяющем уравнении \eqref{lineq}. Для каждого из этих уравнений находились все возможные функции $F$ -- правые части уравнений \eqref{eq}. Только для пяти уравнений из списка (\ref{tst1}--\ref{tst15m}) получен положительный результат. Ниже приводится список найденных уравнений (для них используется спициальная нумерация вида (S...). %, и их дискетные высшие симметрии. 
 Разные уравнения списка (\ref{tst1}--\ref{tst15m}) разделяются с использованием символа ${{\bigstar}}$. %С точностью до замены $x\leftrightarrow y$ список полученных уравнений содержит 4 разныъ

\itemaa 
Для уравнения \eqref{tst5} гиперболическое уравнение имеет вид: 
\mmeq{u_{xy}=2f_a(u_y)u,\quad a\neq 0\label{tzm},} где функция $f_a(u_y)$ удолетворяет уравнению $$(f_a(u_y)+u_y)^2(2f_a(u_y)-u_y)+a^3=0,$$ которое  при $a=1$  совпадает с \eqref{alg1}. Преобразование $$u(x,y)\to u(x,ay),\quad u_y\to au_y,\quad f_a(u_y)\to a{f_1(u_y)},$$ переводит константу $a$ в единицу.

Уравнение \seqref{tzm} было известно \cite{zs01}.
Для уравнения \seqref{tzm} можно показать, что преобразование \cite{ss95}: $$f_a(u_y)+u_y=e^v,\ \ u=2v_x$$ приводит к уравнение Цицейки \eqref{tz} с коэффициентами $$v_{xy}=\frac13e^v-\frac{a^3}6e^{-2v}.$$

\itemaa Для уравнения \eqref{tst5} получаем уравнений Цицейки:
\mmeq{u_{xy}=ae^u+be^{-2u}\label{tzab}.}

\itemaa 
Для уравнения \eqref{tst11} c $\mu=0$ получаем гиперболическое уравнение: 
\mmeq{u_{xy}=2f_a(u_y+b)\sqrt{u_x}\label{tz_neq}.}
В этом случае можно выписать преобразования к предыдущим уравнениям. Замена
 $$u_y=\frac a3{(2e^v+e^{-2v})},\ \ u_x=v_x^2,$$ приводит к уравнению $$ v_{xy}=\frac a3(e^v-e^{-2v});$$ замена
$$u_x=w^2,\ \ u_y=2f_{-a/4^{1/3}}(w_y)+w_y-b,\ \ f_a(u_y+b)=w_y$$ к уравнению $$ w_{xy}=2f_{-a/4^{1/3}}(w_y)w.$$

\itemaa Уравнения \eqref{tst12m} при $\mu=0$ является симметрией двух гиперболических уравнения: 
\mmeq{u_{xy}=2f(u_x)a(u+b),\ \ a\neq 0\label{tzmo},}
\mmeq{u_{xy}=2f(u_x)a\sqrt{u_y+b},\ \ a\neq 0\label{tz_neq1}.}
Заметим, что замены $u(x,y)=\tilde u(x,ay)-b$ и $u(x,y)=\tilde u(x,ay)-by$ приводят к значениям $a=1$ и $b=0$.
При них уравнение \seqref{tzmo} и \seqref{tzm} и уравнения \seqref{tz_neq1} и \seqref{tz_neq} совпадают после преобразования $x\leftrightarrow y$ и вышеуказанных точечных заменах.

\itemaa Для уравнения  \eqref{tst15m} получаем гиперболическое уравнение:
\mmeq{u_{xy}=-2f(u_x)f_a(u_y)\frac{\omega'(u)}{\omega(u)}.\label{tz01}}
Это уравнение появилось в работах \cite{zs01,bzp02}, в них же можно найти преобразование, связывающее уравнение \seqref{tz01} с уравнением Цицейки:
$$e^v=\frac{4c(f(u_x)+u_x)(f_a(u_y)+u_y)\omega(u)}{{\sqrt{c}\omega'(u)-c}}$$
тогда функция $v$ удовлетворяет уравнению $$v_{xy}=e^{v}-4ca^3e^{-2v}.$$ Отметим, что это преобразование в отличие от предыдущих является более сложным -- оно содержит обе частные производные $u_x,u_y$. 

Для нахождения и проверки замен было использовано параметрическое представление кривой третьего порядка \eqref{alg1}, см. \cite{ms12}:
\seq{u_1=\frac13(2e^v+e^{-2v}),\quad f(u_1)=\frac13(e^v-e^{-2v}).} Оно позволяет эффективно вычислять интегралы от рациональных функций $R(u_1,f(u_1))$ по формуле:
\seq{\int R(u_1,f(u_1))du_1=\int R\left(\frac13(2e^v+e^{-2v}),\frac13(e^v-e^{-2v})\right)\frac23(e^v-e^{-2v})dv}

Высшие симметрии уравнений \seqref{tzm}-\seqref{tz01} в $y$ направлении также являются уравнениями пятого порядка. Уравнения \seqref{tzab} и \seqref{tz01} симметричны при замене $x\leftrightarrow y$
\qed

В теореме получено шесть уравнений, однако две их пары эквивалентны друг другу с точностью до преобразования $x\leftrightarrow y$ и простых точечных замен. Приведем список четырех уравнений в наиболее простой форме, где исключены константы:
\seqs{&u_{xy}=e^u+e^{-2u},\\ &u_{xy}=2f(u_x)u, \\ &u_{xy}=2f(u_x)\sqrt{u_y},\\ &u_{xy}=-2f(u_x)f(u_y)\frac{\omega'(u)}{\omega(u)}.}

\end{document}